\documentclass{article}
\usepackage{amsmath}
\usepackage{amssymb}
\usepackage{amsthm}
\allowdisplaybreaks
\begin{document}
\title{\bf On the Hybrid Mean Value of
 Generalized Dedekind Sums, Generalized  Hardy
Sums and Kloosterman Sums \footnote{This work is supported by the
Natural Science Foundation of the Education Department of Shannxi
Province of China (No.16JK1456) }}
\author{ Qing Tian \footnote{E-mail: qingtian@xauat.edu.cn( Qing Tian)} \\{\small
School of Science, Xi'an University of Architecture and Technology }\\{ \small Xi'an, 710055, Shaanxi, P. R. China }\\
}
\date{}
\maketitle \baselineskip 16pt \begin{center}
\begin{minipage}{120mm}
{\bf Abstract}\\The main purpose of this paper is to study the
hybrid mean value problem involving generalized Dedekind sums,
generalized Hardy sums and Kloosterman sums, and give some exact
computational formulae for them by using the properties of Gauss
sums and the mean value theorem of the Dirichlet L-function.
\ \ {\small  }\\
 {\bf Keywords:}\ \ Hybrid mean value, Kloosterman sums, Generalized Dedekind sums, Generalized Hardy sums.\\
{\bf MSC: } 11F20, 11L05
\end{minipage}
\end{center}

\section*{ 1. Introduction }
  Suppose that $k$ is a positive integer, then for an arbitrary
integer $h,m,n$, The generalized  Dedekind sums are defined by
$$S(h,m,n,k)=\sum^{k}_{j=1}\overline{B}_{m}\left(\frac{j}{k}\right)\overline{B}_{n}\left(\frac{hj}{k}\right),$$
where
$$\overline{B}_{m}(x)=\left\{
                    \begin{array}{ll}
                       B_{m}(x-[x]), & \hbox{if \textit{x} is not an integer;} \\
                      0, & \hbox{if \textit{x} is an integer.}
                    \end{array}
                  \right.
$$
${B}_{m}(x)$ is the Bernoulli polynomial, $\overline{B}_{m}(x)$
defined on the interval $0< x\leq 1$ is the $n$-th Bernoulli
periodic function. For $m=n=1$, $S(h,1,1,q)=S(h,q)$ are the
classical Dedekind sums, which were studied by many experts because
of the prominent role they play in the transformation theory of the
Dedekind eta-function. Some arithmetical properties of $S(h,q)$ can
be found in Apostol [1] and Carlitz [2]. The most famous property of
the Dedekind sums may be the reciprocity formula ([3][4])
$$S(h,k)+S(k,h)=\frac{h^2+k^2+1}{12hk}-\frac{1}{4}.$$

In [5], Berndt gave certain sums called Hardy sums which are
related to the Dedekind sums, and also obtained some arithmetic
properties (see [6]). Sitaramachandrarao [7] and Pettet [8] used
elementary methods to express the Hardy sums in terms of the
Dedekind sums. H.N.Liu [9] generalize the Hardy sums as
follows
\begin{eqnarray*}
& &s_{1}(h,m,k)=\sum^{k}_{j=1}(-1)^{[\frac{hj}{k}]}\overline{B}_{m}\left(\frac{j}{k}\right),\quad \quad s_{2}(h,m,n,k)=\sum^{k}_{j=1}(-1)^j\overline{B}_{m}\left(\frac{j}{k}\right)\overline{B}_{n}\left(\frac{hj}{k}\right),\\
&
&s_{3}(h,n,k)=\sum^{k}_{j=1}(-1)^{j}\overline{B}_{n}\left(\frac{hj}{k}\right),\quad
\quad \quad
s_{5}(h,m,k)=\sum^{k}_{j=1}(-1)^{j+[\frac{hj}{k}]}\overline{B}_{m}\left(\frac{j}{k}\right).
\end{eqnarray*}
For $m=n=1$, the sums $s_{1}(h,k)=s_{1}(h,1,k), s_{2}(h,k)=s_{2}(h,1,1,k), s_{3}(h,k)=s_{3}(h,1,k)$ and $s_{5}(h,k)=s_{5}(h,1,k)$ are classical Hardy sums defined in [5]. H. N. Liu's research paper [9] also express the generalized Hardy sums in term of generalized
Dedekind sums, that is

{\bf Proposition 1.1} Let $h$, $q$ be positive integer with $(h,q)=1$, then
$$\begin{array}{ll}
s_{1}(h,m,q)=2\cdot S(h,m,1,q)-4\cdot S\left(\frac{h}{2},m,1,q\right),& \hbox{if $h$ is even number} \\
s_{2}(h,m,n,q)=2^m\cdot S(2h,m,n,q)-S(h,m,n,q), & \hbox{if $q$ is even number}\\
s_{3}(h,n,q)=2\cdot S(h,,1,n,q)-4\cdot S(2h,1,n,q), & \hbox{if $q,n$ are odd number}\\
s_{5}(h,m,q)=2^{m+1}\cdot S(2h,m,1,q)+2^{m+1}\cdot S(h,m,1,2q)
\\\quad \quad \quad \quad \quad \quad  -(2+2^{m+2})\cdot S(h,m,1,q), & \hbox{if $h+q$ is even number}\\
\end{array}$$
where $\bar{2}\cdot 2\equiv 1\bmod q$. Moreover, each one of
$$\left\{
                    \begin{array}{ll}
                      s_1(h,m,q)\quad (h+m \quad  \textrm{even}),\quad\quad s_2(h,m,n,q)\quad (h+m+q \quad  \textrm{odd}) \\
                      s_3(h,n,q)\quad (h+q \quad  \textrm{odd}), \quad\quad s_5(h,m,k)\quad (h+m+q \quad \textrm{even}).
                    \end{array}
                  \right.
$$
is zero.

Recently, some authors studied the hybrid mean value of Dedekind sums or Hardy sums with Kloosterman sums defined by
$$K(n,q)=\mathop{{\sum}'}^q_{c=1}e\left(\frac{mc+\bar{c}}{q}\right),$$
where $\mathop{{\sum}'}^q_{c=1}$ denotes the summation over all $c$ such that $(c,q)=1$, $e(y)=\exp(2\pi iy)$ and $\bar{c}\cdot c\equiv 1\bmod q$. And they found there are some close
relationships between the functions. Y. N. Liu et al. [10] gave several explicit formulae for
$$\sum_{a=1}^{p-1}\sum_{b=1}^{p-1}K^2(a,p) K^2(b,p)S^{k}(a\bar{b},p)$$ under the condition $q=p$ is a prime.
H. Zhang et al. [11] and W. Peng et al. [12] also obtained identities for
$$\sum_{a=1}^{p-1}\sum_{b=1}^{p-1}K(a,p) K(b,p)s_1(2a\bar{b},p).$$ and $$\sum_{a=1}^{p-1}\sum_{b=1}^{p-1}K(a,p) K(b,p)s_5(a\bar{b},p)$$ respectively.

Naturally, one might consider whether the hybrid mean value be extended to generalized Dedekind sums
$S(h,m,n,q)$ or certain generalized Hardy sums with Kloosterman
sums $K(n,q)$ under the condition of composite number $q$? If yes, then what can be expected? These
problems may be interesting. In this paper, we shall study the problems  and give some exact
computational formulae by using the prosperities of Gauss sums and the mean value theorem of the Dirichlet L-function. That is, we shall prove the following:

{\bf Theorem 1.} Let $q$ be a square-full number, $m\equiv n\equiv
1\bmod 2$. Then we have
\begin{eqnarray*}
&&\mathop{{\sum}'}_{a=1}^{q}\mathop{{\sum}'}_{b=1}^{q}K(a,q) K(b,q) S(\bar{a}b,m,n,q)\\
&&=q^{4-2m-2n}\cdot \sum^{m+n}_{l=0}q^{l}\cdot {r_{m, n, l}\cdot
\prod_{p|q}\left(1-\frac{1}{p}\right)\left(1-\frac{1}{p^l}\right)\left(1-\frac{1}{p^{l-m-n+1}}\right)}
\end{eqnarray*}
where
$$r_{m,n,l}=B_{m+n-l}\displaystyle\mathop{\sum^{m}_{a=0}\sum^{n}_{b=0}}_{a+b\geq
m+n-l}B_{m-a}B_{n-b}\frac{\left(\begin{array}{c}
                                                                                                                        m \\
                                                                                                                         a
                                                                                                                       \end{array}
\right)\left(\begin{array}{c}
                      n \\
                      b
                    \end{array}\right)
\left(\begin{array}{c}
        a+b+1 \\
        m+n-l
      \end{array}
\right)}{a+b+1},$$ $B_{m}$ is Bernoulli number,  $\left(
                                    \begin{array}{c}
                                      m \\
                                      a \\
                                    \end{array}
                                  \right)
=\frac{m!}{a!(m-a)!}$, $\phi_{l}(q)=\prod_{p|q}(1-\frac{1}{p^l})$, $\prod_{p|q}$ denotes the products of all prime divisors of $q$
and $\phi(q)=q\phi_1(q).$

{\bf Theorem 2.}  Let $q$ be a square-full number, $m\equiv 1\bmod
2$. Then we have
\begin{eqnarray*}
&&\mathop{{\sum}'}_{a=1}^{q}\mathop{{\sum}'}_{b=1}^{q}K(a,q)K(b,q)s_1(2\bar{a}b,m,q)\\
&&=q^{m-2}\cdot \sum^{m+1}_{l=0}q^{l-m}\cdot r_{m,1,l}\cdot
\frac{-2^l-2}{2^{m-1}+1}\cdot
\prod_{p|q}\left(1-\frac{1}{p}\right)\left(1-\frac{1}{p^l}\right)\left(1-\frac{1}{p^{l-m}}\right)
\end{eqnarray*}

{\bf Theorem 3.}  Let $q$ be a square-full even number, $m\equiv
n\equiv 1\bmod 2$. Then we have
\begin{eqnarray*}
&&\mathop{{\sum}'}_{a=1}^{q}\mathop{{\sum}'}_{b=1}^{q}K(a,q)K(b,q)s_2(\bar{a}b,m,n,q)\\
&&=q^{2-2m}\cdot \sum^{m+1}_{l=0}r_{m,1,l}\cdot (2^m\frac{2^m-2^{l-1}+1}{2^{m-1}+1}-1)\cdot
q^{l-m}\cdot\prod_{p|q}\left(1-\frac{1}{p}\right)\left(1-\frac{1}{p^l}\right)\left(1-\frac{1}{p^{l-m}}\right)
\end{eqnarray*}

{\bf Theorem 4.} Let $q$ be a square-full odd number, $n\equiv
1\bmod 2$. Then we have
\begin{eqnarray*}
&&\mathop{{\sum}'}_{a=1}^{q}\mathop{{\sum}'}_{b=1}^{q}K(a,q)K(b,q)s_3(\bar{a}b,n,q)\\
&&=q^{2-2n}\cdot\sum^{n+1}_{l=0}r_{1,n,l}\cdot q^{l}\left(2-4\cdot
\frac{2^l-2^{n+1}-2}{2^{n}+2}\right)\cdot
\prod_{p|q}\left(1-\frac{1}{p}\right)\left(1-\frac{1}{p^l}\right)\left(1-\frac{1}{p^{l-n}}\right)
\end{eqnarray*}

{\bf Theorem 5.} Let $q$ be a square-full odd number, $m\equiv
1\bmod 2$. Then we have
\begin{eqnarray*}
&&\mathop{\mathop{{\sum}'}_{a=1}^{q}}_{(2a-1,q)=1}\mathop{\mathop{{\sum}'}_{b=1}^{q}}_{(2b-1,q)=1}K(2a-1,q)K(2b-1,q)s_5((\overline{2a-1})(2b-1),m,q)\\
&&=q^{2-2m}\cdot \sum^{m+1}_{l=0}r_{m,1,l}\cdot q^{l}\left(2-4\cdot
\frac{2^l-2^{m+1}-2}{2^{m}+2}\right)
\cdot \prod_{p|q}\left(1-\frac{1}{p}\right)\left(1-\frac{1}{p^l}\right)\left(1-\frac{1}{p^{l-m}}\right).
\end{eqnarray*}

The present work is a generalization of [11] and [12].

For general number $q>2 $ , we can only get some asymptotic
formulae, whether there exits  the identities for the hybrid mean
value for these sums are open problems.

\section*{ 2. Several Lemmas}

Before starting our proof of the theorems, several lemmas will be useful.

{\bf Lemma 2.1.} Let $h$, $q$ be positive integer with $q\geq 3$ and $(h,q)=1$, $m\equiv n\equiv 1\bmod 2$. Then
we have
$$S(h,m,n,q)=\frac{-4m!n!}{(2\pi i)^{m+n}q^{m+n-1}}\cdot \sum_{d|q}\frac{d^{m+n}}{\phi(d)}\mathop{\sum_{\chi\bmod d}}_{\chi(-1)=-1}\overline{\chi}(h)L(m,\chi)L(n,\overline{\chi}),$$
 where$\sum_{d|q}$ denotes the sums over all divisors of $q$ and $L(m,\chi)$ denotes the Dirichlet L-function corresponding to character $\chi\bmod d$.

{\bfseries Proof.} See Theorem 2.3 of [9].

{\bf Lemma 2.2.}  Let $q\geq 3$ be an odd number. Then for odd numbers
$h,q$ with $(h,q)=1$, we have
\begin{eqnarray*}
s_{5}(h,m,q)=2\cdot S(h,m,1,q)-4\cdot S(\bar{2}h,m,1,q).
\end{eqnarray*}

{\bfseries Proof.} From the Proposition 1.1, we know that if $h+q$
is even number, generalized Hardy sums $s_5(h,m,q)$ can be expressed in term of generalized Dedekind sums, that is
\begin{eqnarray}
s_{5}(h,m,q)=2^{m+1}\cdot S(2h,m,1,q)+2^{m+1}\cdot
S(h,m,1,2q)-(2+2^{m+2})\cdot S(h,m,1,q).
\end{eqnarray}
Now we simplify the formula (1). Considering the second part $S(h,m,1,2q)$ firstly, by using Lemma 1 we have
\begin{eqnarray*}
&&S(h,m,1,2q)\\
&&=-\frac{4m!}{(2\pi
i)^{m+1}(2q)^{m}}\cdot \sum_{d|2q}\frac{d^{m+1}}{\phi(d)}\mathop{\sum_{\chi\bmod
d}}_{\chi(-1)=-1}\overline{\chi}(h)L(m,\chi)L(1,\overline{\chi})\\
&&=-\frac{4m!}{(2\pi
i)^{m+1}(2q)^{m}}\cdot \left(\sum_{d|q}\frac{(2d)^{m+1}}{\phi(2d)}\mathop{\sum_{\chi\bmod
2d}}_{\chi(-1)=-1}\overline{\chi}(h)L(m,\chi)L(1,\overline{\chi})\right.\\
&&\qquad
\left.+\sum_{d|q}\frac{d^{m+1}}{\phi(d)}\mathop{\sum_{\chi\bmod
d}}_{\chi(-1)=-1}\overline{\chi}(h)L(m,\chi)L(1,\overline{\chi})\right)\\
&&=-\frac{4m!}{(2\pi i)^{m+1}(2q)^{m}}\cdot \left(2^{m+1}\cdot
\sum_{d|q}\frac{d^{m+1}}{\phi(d)}\mathop{\sum_{\chi\bmod
d}}_{\chi(-1)=-1}\overline{\chi}\chi_{2}^{0}(h)L(m,\chi\chi_{2}^{0})L(1,\overline{\chi}\chi_{2}^{0})\right.\\
&&\left.+\sum_{d|q}\frac{d^{m+1}}{\phi(d)}\mathop{\sum_{\chi\bmod
d}}_{\chi(-1)=-1}\overline{\chi}(h)L(m,\chi)L(1,\overline{\chi})\right)\\
\end{eqnarray*}
where $\chi_{2}^{0}$ denotes the principal character modulo 2.

From the Euler infinite product formula (see Theorem 11.6 of [13]),
we have
\begin{eqnarray*}
L(m,\chi\chi_{2}^{0})&=&\prod_{p_1}\left(1-\frac{\chi(p_1)\chi_{2}^{0}(p_1)}{p_1^{m}}\right)^{-1}
=\prod_{p_1>2}\left(1-\frac{\chi(p_1)}{p_1^{m}}\right)^{-1}\\
&=&\left(1-\frac{\chi(2)}{2^m}\right)\prod_{p_1}\left(1-\frac{\chi(p_1)}{p_1^{m}}\right)^{-1}=\left(1-\frac{\chi(2)}{2^m}\right)L(m,\chi)
\end{eqnarray*}
\begin{eqnarray*}
L(1,\overline{\chi}\chi_{2}^{0})=\prod_{p_2}\left(1-\frac{\overline{\chi}(p_2)\chi_{2}^{0}(p_2)}{p_2}\right)^{-1}
=\left(1-\frac{\overline{\chi}(2)}{2}\right)L(1,\overline{\chi})
\end{eqnarray*}
where $\displaystyle\prod_p$ denotes the product over all primes $p$.

That is we have the identity
\begin{eqnarray}
&&S(h,m,1,2q)\nonumber\\
&&=-\frac{4m!}{(2\pi i)^{m+1}q^{m}}\cdot \frac{1}{2^m}\cdot
\left(2^{m+1}\cdot\sum_{d|q}\frac{d^{m+1}}{\phi(d)}\mathop{\sum_{\chi\bmod
d}}_{\chi(-1)=-1}\left(1-\frac{\chi(2)}{2^m}\right)\left(1-\frac{\overline{\chi}(2)}{2}\right)\overline{\chi}(h)L(m,\chi)L(1,\overline{\chi})\right.\nonumber\\
&&\qquad\left.+\sum_{d|q}\frac{d^{m+1}}{\phi(d)}\mathop{\sum_{\chi\bmod
d}}_{\chi(-1)=-1}\bar{\chi}(h)L(m,\chi)L(1,\bar{\chi})\right)\nonumber\\
&&=-\frac{4m!}{(2\pi
i)^{m+1}q^{m}}\cdot\sum_{d|q}\frac{d^{m+1}}{\phi(d)}\mathop{\sum_{\chi\bmod
d}}_{\chi(-1)=-1}\left(2+\frac{1}{2^{m-1}}-\overline{\chi}(2)-\frac{\chi(2)}{2^{m-1}}\right)\overline{\chi}(h)L(m,\chi)L(1,\overline{\chi})\nonumber\\
&&=(2+\frac{1}{2^{m-1}})\cdot S(h,m,1,q)-S(2h,m,1,q)-\frac{1}{2^{m-1}}\cdot S(\bar{2}h,m,1,q)
\end{eqnarray}

Combining (1) with (2), it follows that
\begin{eqnarray*}
s_{5}(h,m,q)=2\cdot S(h,m,1,q)-4\cdot S(\bar{2}h,m,1,q).
\end{eqnarray*}
This proves Lemma 2.2.

 {\bf Lemma 2.3.} Let $q\geq 3$ be an integer and $\chi$ be a non-principal character mod
 $q$. Then we have
 $$\mathop{{\sum}'}_{a=1}^{q}\chi(a)K(a,q)=\tau^2(\chi)$$
{\bfseries Proof.} From the properties of reduced
residue system, it is known that if $a$ pass through a reduced
residue system mod $q$, then for any integer $c$ with $(c,q)=1$, $ac
$ also pass through a reduced residue system mod $q$, by the  definition of Gauss sums, we have
\begin{eqnarray*}
\mathop{{\sum}'}_{a=1}^{q}\chi(a)K(a,q)&=&\mathop{{\sum}'}_{a=1}^{q}\chi(a)\mathop{{\sum}'}_{c=1}^{q}e\left(\frac{ac+\bar{c}}{q}\right)\\
&=&\mathop{{\sum}'}_{c=1}^{q}e\left(\frac{\bar{c}}{q}\right)\mathop{{\sum}'}_{c=1}^{q}\chi(a)e\left(\frac{ac}{q}\right)\\
&=&\mathop{{\sum}'}_{c=1}^{q}\chi(\bar{c})e\left(\frac{\bar{c}}{q}\right)\mathop{{\sum}'}_{a=1}^{q}\chi(ac)e\left(\frac{ac}{q}\right)\\
&=&\tau^2(\chi).
\end{eqnarray*}
This proves Lemma 2.3.

{\bf Lemma 2.4.} Let $q$ be a square-full number. Then for any
non-primitive character $\chi\bmod q$, we have the identity
$$\tau(\chi)=\sum_{a=1}^{q}\chi(a)e\left(\frac{a}{q}\right)=0$$

{\bfseries Proof.} It is known that $\tau^2(\chi)$ is a
multiplicative function, so without loss of generality we  assume
that $q=p^\alpha$, where $p$ is a prime and $\alpha\geq 2$. If
$\chi$ is a non-primitive character modulo $q=p^{\alpha}$, then
$\chi$ must be a character modulo  $p^{\alpha-1}$. Note that the
trigonometric identity $\sum_{a=0}^{p-1}e(\frac{a}{p})=0$. From the
properties of the reduced residue system modulo $p^{\alpha-1}$, it is easy to get
\begin{eqnarray*}
\tau(\chi)&=&\sum_{a=1}^{q}\chi(a)e\left(\frac{a}{q}\right)=\sum_{a=1}^{p-1}\sum_{b=1}^{p^{\alpha-1}}\chi(b)e\left(\frac{ap^{\alpha-1}+b}{p^{\alpha}}\right)\\
&=&\sum_{b=1}^{p^{\alpha-1}}\chi(b)e\left(\frac{b}{p^\alpha}\right)\sum_{a=1}^{p-1}e\left(\frac{a}{p}\right)=0.
\end{eqnarray*}
This proves Lemma 2.4.

{\bf Lemma 2.5.} Let $q\geq 2$ be an integer, $m\equiv n\equiv 1$. We
have
$$\mathop{\sum_{\chi\bmod
q}}_{\chi(-1)=-1}L(m,\chi)L(n,\overline{\chi})=-\frac{(2\pi
i)^{m+n}\phi(q)}{4m!n!}\cdot \left(\sum^{m+n}_{l=0}r_{m,n,l}\cdot
\phi_{l}(q)\cdot q^{l-m-n}-\frac{B_mB_n\phi_{m+n-1}(q)}{q}\right).$$

{\bfseries Proof.}  See Theorem 3 of [9].

{\bf Lemma 2.6.} Let $q$ be square-full number. Then we have
\begin{eqnarray*}
&&\displaystyle\mathop{\mathop{{{\sum}^*}}_{\chi\bmod q}}_{\chi(-1)=-1}L(m,\chi)L(n,\overline{\chi})\\
&&=-\frac{(2\pi i
)^{m+n}}{4m!n!}\sum^{m+n}_{l=0}{r_{m,n,l}\cdot q^{l-m-n+1}\cdot \prod_{p|q}\left(1-\frac{1}{p}\right)\left(1-\frac{1}{p^l}\right)\left(1-\frac{1}{p^{l-m-n+1}}\right)},
\end{eqnarray*}
where $\displaystyle\mathop{\mathop{{{\sum}^*}}_{\chi\bmod q}}_{\chi(-1)=-1}$ denotes the sums over all odd primitive characters mod $q$.

{\bfseries Proof.} Noting that $q$ is square-full number and
\begin{eqnarray*}\displaystyle\mathop{\mathop{{{\sum}}}_{\chi\bmod q}}_{\chi(-1)=-1}L(m,\chi)L(n,\overline{\chi})=\sum_{d|q}\displaystyle\mathop{\mathop{{{\sum}^*}}_{\chi\bmod
d}}_{\chi(-1)=-1}L(m,\chi\chi^{0}_{q})L(n,\overline{\chi}\chi^{0}_{q}),
\end{eqnarray*} by using M\"{o}bius inverse formula, we have
\begin{eqnarray*}
\displaystyle\mathop{\mathop{{{\sum}^*}}_{\chi\bmod
q}}_{\chi(-1)=-1}L(m,\chi)L(n,\overline{\chi})&=&\displaystyle\mathop{\mathop{{{\sum}^*}}_{\chi\bmod q}}_{\chi(-1)=-1}L(m,\chi\chi^{0}_{q})L(n,\overline{\chi}\chi^{0}_{q})\\
&=&\sum_{d|q}\mu(d)\displaystyle\mathop{\mathop{{{\sum}}}_{\chi\bmod
{\frac{q}{d}}}}_{\chi(-1)=-1}L(m,\chi)L(n,\overline{\chi}).
\end{eqnarray*}
According to Lemma 2.5, we get
\begin{eqnarray*}
\displaystyle\mathop{\mathop{{{\sum}^*}}_{\chi\bmod
q}}_{\chi(-1)=-1}L(m,\chi)L(n,\overline{\chi})&=&-\frac{(2\pi
i)^{m+n}}{4m!n!}\cdot \left[\sum^{m+n}_{l=0}r_{m,n,l}\cdot \sum_{d|q}\mu(d)\phi(\frac{q}{d})\phi_{l}\left(\frac{q}{d}\right)\left(\frac{q}{d}\right)^{l-m-n}\right.\\
&&\left.+B_mB_n\cdot \sum_{d|q}\mu(d)\phi(\frac{q}{d})\frac{\phi_{m+n-1}(\frac{q}{d})}{\frac{q}{d}}\right].
\end{eqnarray*}
Since $\mu(n)$, $\phi(n)$ and
$\phi_{l}(n)$ are multiplicative number, we have
\begin{eqnarray*}
&&\sum_{d|q}\mu(d)\phi\left(\frac{q}{d}\right)\phi_{l}\left(\frac{q}{d}\right)\left(\frac{q}{d}\right)^{l-m-n}
=\prod_{p|q}\sum_{d|p^\alpha}\mu(d)\phi\left(\frac{p^\alpha}{d}\right)\phi_{l}\left(\frac{p^\alpha}{d}\right)\left(\frac{p^\alpha}{d}\right)^{l-m-n}\\
&&=q^{l-m-n+1}\cdot\prod_{p|q}\left(1-\frac{1}{p}\right)\left(1-\frac{1}{p^l}\right)\left(1-\frac{1}{p^{l-m-n+1}}\right)
.\end{eqnarray*} Using the same methods, we get
$$\sum_{d|q}\mu(d)\phi\left(\frac{q}{d}\right)\phi_{m+n-1}\left(\frac{q}{d}\right)\frac{d}{q}=0.$$
Due to the  discussion above, we obtain
\begin{eqnarray*}
&&\displaystyle\mathop{\mathop{{{\sum}^*}}_{\chi\bmod q}}_{\chi(-1)=-1}L(m,\chi)L(n,\overline{\chi})\\
&&=-\frac{(2\pi i
)^{m+n}}{4m!n!}\cdot\sum^{m+n}_{l=0}{r_{m,n,l}\cdot q^{l-m-n+1}\cdot\prod_{p|q}\left(1-\frac{1}{p}\right)\left(1-\frac{1}{p^l}\right)\left(1-\frac{1}{p^{l-m-n+1}}\right)}
.\end{eqnarray*}
This proves  Lemma 2.6.

 {\bf Lemma 2.7.} Let $q$ be square-full number, then
we have the identity
\begin{eqnarray*}
&&\displaystyle\mathop{\mathop{{{\sum}^*}}_{\chi\bmod q}}_{\chi(-1)=-1}L(m,\chi\chi^{0}_{2})L(n,\overline{\chi}\chi^{0}_{2})\\
&&=-\frac{(2\pi
i)^{m+n}}{4m!n!}\cdot \sum^{m+n}_{l=0}r_{m,n,l}\cdot \left(1-\frac{1}{2^l}\right)\cdot 2^{l-m-n}\cdot q^{l-m-n+1}\cdot
\prod_{p|q}\left(1-\frac{1}{p}\right)\left(1-\frac{1}{p^l}\right)\left(1-\frac{1}{p^{l-m-n+1}}\right).
\end{eqnarray*}
{\bfseries Proof.} Note that
\begin{eqnarray*}
\mathop{\sum_{\chi\bmod 2q}}_{\chi(-1)=-1}L(m,\chi)L(n,\overline{\chi})&=&\mathop{\sum_{\chi\bmod q}}_{\chi(-1)=-1}L(m,\chi\chi^{0}_{2})L(n,\overline{\chi}\chi^{0}_{2})\nonumber\\
&=&\sum_{d|q}\displaystyle\mathop{\mathop{{{\sum}^*}}_{\chi\bmod
d}}_{\chi(-1)=-1}L(m,\chi\chi^{0}_{2q})L(n,\overline{\chi}\chi^{0}_{2q})
,\end{eqnarray*} and  M\"{o}bius inverse formula, we have
\begin{eqnarray*}
\displaystyle\mathop{\mathop{{{\sum}^*}}_{\chi\bmod q}}_{\chi(-1)=-1}L(m,\chi\chi^{0}_{2})L(n,\overline{\chi}\chi^{0}_{2})&=&\displaystyle\mathop{\mathop{{{\sum}^*}}_{\chi\bmod q}}_{\chi(-1)=-1}L(m,\chi\chi^{0}_{2q})L(n,\overline{\chi}\chi^{0}_{2q})\\
&=&\sum_{d|q}\mu(d)\mathop{\sum_{\chi\bmod \frac{q}{d}}}_{\chi(-1)=-1}L(m,\chi\chi^{0}_{2q})L(n,\overline{\chi}\chi^{0}_{2q})\\
&=&\sum_{d|q}\mu(d)\mathop{\sum_{\chi\bmod\frac{q}{d}}}_{\chi(-1)=-1}L(m,\chi\chi^{0}_{\frac{2q}{d}})L(n,\overline{\chi}\chi^{0}_{\frac{2q}{d}})\nonumber\\
&=&\sum_{d|q}\mu(d)\mathop{\sum_{\chi\bmod\frac{2q}{d}}}_{\chi(-1)=-1}L(m,\chi)L(n,\overline{\chi}).
\end{eqnarray*}

According to the Lemma 2.6, it follows that
\begin{eqnarray*}
&&\displaystyle\mathop{\mathop{{{\sum}^*}}_{\chi\bmod  q}}_{\chi(-1)=-1}L(m,\chi\chi^{0}_{2})L(n,\overline{\chi}\chi^{0}_{2})\\
&&=-\frac{(2\pi
i)^{m+n}}{4m!n!}\cdot \left[\sum^{m+n}_{l=0}r_{m,n,l}\cdot \sum_{d|q}\mu(d)\phi\left(\frac{2q}{d}\right)\phi_{l}\left(\frac{2q}{d}\right)\left(\frac{2q}{d}\right)^{l-m-n+1}\right.\\
&&\left.\quad +\sum_{d|q}\mu(d)\phi\left(\frac{2q}{d}\right)\cdot \frac{B_mB_n\phi_{m+n-1}\left(\frac{2q}{d}\right)}{\frac{2q}{d}}\right]\\
&&=-\frac{(2\pi
i)^{m+n}}{4m!n!}\cdot \sum^{m+n}_{l=0}r_{m,n,l}\cdot \left(1-\frac{1}{2^l}\right)\cdot 2^{l-m-n}\cdot q^{l-m-n+1}\cdot
\prod_{p|q}\left(1-\frac{1}{p}\right)\left(1-\frac{1}{p^l}\right)\left(1-\frac{1}{p^{l-m-n+1}}\right).\\
\end{eqnarray*}
This proves Lemma 2.7.

{\bf Lemma 2.8.} Let $q$ be square-full odd number, we have
\begin{eqnarray*}
&&\displaystyle\mathop{\mathop{{{\sum}^*}}_{\chi\bmod q}}_{\chi(-1)=-1}\overline{\chi}(2)L(m,\chi)L(n,\overline{\chi})\\
&&=-\frac{(2\pi
i)^{m+n}}{4m!n!}\cdot \sum^{m+n}_{l=0}r_{m,n,l}\cdot q^{l-m-n+1}\cdot \frac{2^l-2^{m+n}-2}{2^m+2^n}
\cdot \prod_{p|q}\left(1-\frac{1}{p}\right)\left(1-\frac{1}{p^l}\right)\left(1-\frac{1}{p^{l-m-n+1}}\right).
\end{eqnarray*}

{\bfseries Proof.} From the proof of Lemma 2.2, we know that
\begin{eqnarray*}
L(m,\chi\chi^{0}_{2})L(n,\overline{\chi}\chi^{0}_{2})
&=&L(m,\chi)L(n,\overline{\chi})\left(1-\frac{\chi(2)}{2^m}\right)\left(1-\frac{\overline{\chi}(2)}{2^n}\right)\\
&=&L(m,\chi)L(n,\overline{\chi})\left[1+\frac{1}{2^{m+n}}-\left(\frac{\chi(2)}{2^m}+\frac{\overline{\chi}(2)}{2^n}\right)\right].
\end{eqnarray*} Note that
$\displaystyle\mathop{\mathop{{\sum}^*}_{\chi \bmod
q}}_{\chi(-1)=-1}\chi(2)=\displaystyle\mathop{\mathop{{\sum}^*}_{\chi
\bmod q}}_{\chi(-1)=-1}\overline{\chi}(2)$ we get
\begin{eqnarray*}
&&\displaystyle\mathop{\mathop{{{\sum}^*}}_{\chi\bmod q}}_{\chi(-1)=-1}\overline{\chi}(2)L(m,\chi)L(n,\overline{\chi})\\
&&=\displaystyle\mathop{\mathop{{{\sum}^*}}_{\chi\bmod q}}_{\chi(-1)=-1}\left(\frac{1}{2^n}+\frac{1}{2^m}\right)^{-1}\cdot \left(L(m,\chi\chi^{0}_{2})L(n,\overline{\chi}\chi^{0}_{2})-\left(1+\frac{1}{2^{m+n}}\right)\cdot L(m,\chi)L(n,\overline{\chi})\right)\\
&&=-\frac{(2\pi
i)^{m+n}}{4m!n!}\cdot \sum^{m+n}_{l=0}r_{m,n,l}\cdot q^{l-m-n+1}\cdot \frac{2^{l-1}-2^{m+n-1}-1}{2^{m-1}+2^{n-1}}\cdot
\prod_{p|q}\left(1-\frac{1}{p}\right)\left(1-\frac{1}{p^l}\right)\left(1-\frac{1}{p^{l-m-n+1}}\right)
.\end{eqnarray*}
This proves Lemma 2.8.

\section*{\S 3. Proof of the theorems}

In this section, we shall complete the proof of the theorems.

First we give a hybrid mean value formula for generalized Dedekind
sums with Kloosterman sums. Note that if $\chi$ is primitive character mod $q$, the Gauss sums $\tau(\chi)=\sqrt{q}$ and
$$\left|\mathop{{\sum}'}_{a=1}^{q}\chi(a)K(a,q)\right|=|\tau^2(\chi)|=q.$$
From Lemma 2.3 and 2.4 we have if $q$ is a square-full  number and $m\equiv n\equiv 1(\bmod 2)$
\begin{eqnarray*}
&&\mathop{{\sum}'}_{a=1}^{q}\mathop{{\sum}'}_{b=1}^{q}K(a,q)K(b,q)S(\bar{a}b,m,n,q)\\
&&=-\frac{4m!n!}{(2\pi
i)^{m+n}q^{m+n-1}}\cdot \sum_{d|q}\frac{d^{m+n}}{\phi(d)}\mathop{\sum_{\chi\bmod
d}}_{\chi(-1)=-1}\mathop{{\sum}'}_{a=1}^{q}\mathop{{\sum}'}_{a=1}^{q}K(a,q)K(b,q)\bar{\chi}(\bar{a}b)L(m,\chi)L(n,\bar{\chi})\\
&&=-\frac{4m!n!}{(2\pi
i)^{m+n}q^{m+n-1}}\cdot \sum_{d|q}\frac{d^{m+n}}{\phi(d)}\mathop{\sum_{\chi\bmod
d}}_{\chi(-1)=-1}\left|\mathop{{\sum}'}_{a=1}^{q}\chi(a)K(a,q)\right|^2\cdot L(m,\chi)L(n,\bar{\chi})\\
&&=-\frac{4m!n!}{(2\pi
i)^{m+n}q^{m+n-1}}\cdot \frac{q^{m+n}}{\phi(q)}\cdot\mathop{\mathop{{\sum}^*}_{\chi\bmod
q}}_{\chi(-1)=-1}q^2\cdot L(m,\chi)L(n,\bar{\chi})\\
&&=\frac{q^{4-m-n}}{\phi(q)}\cdot \sum^{m+n}_{l=0}{r_{m,n,l}\cdot q^{l}\cdot \prod_{p|q}\left(1-\frac{1}{p}\right)\left(1-\frac{1}{p^l}\right)\left(1-\frac{1}{p^{l-m-n+1}}\right)}.
\end{eqnarray*}
This completes the proof of Theorem 1.

Now we give some hybrid mean value formulae for generalized Hardy sums and Kloosterman sums. That is we will prove Theorem 2-5.

From the proposition 1.1, Lemma 2.6 and Lemma 2.8 together, we have if $q$ is a square-full  number and $m\equiv 1(\bmod 2)$
\begin{eqnarray*}
&&\mathop{{\sum}'}_{a=1}^{q}\mathop{{\sum}'}_{b=1}^{q}K(a,q)K(b,q)s_1(2\bar{a}b,m,q)\\
&&=2\mathop{{\sum}'}_{a=1}^{q}\mathop{{\sum}'}_{b=1}^{q}K(a,q)K(b,q)S(2\bar{a}b,m,1,q)-4\mathop{{\sum}'}_{a=1}^{q}\mathop{{\sum}'}_{b=1}^{q}K(a,q)K(b,q)S(\bar{a}b,m,1,q)\\
&&=-\frac{4m!}{(2\pi
i)^{m+1}q^{m}}\cdot \sum_{d|q}\frac{d^{m+1}}{\phi(d)}\mathop{\sum_{\chi\bmod
d}}_{\chi(-1)=-1}\mathop{{\sum}'}_{a=1}^{q}\left|\chi(a)K(a,q)\right|^2\cdot(2\bar{\chi}(2)-4)\cdot L(m,\chi)L(1,\bar{\chi})\\
&&=-\frac{4m!}{(2\pi
i)^{m+1}q^{m}}\cdot \frac{q^{m+1}}{\phi(q)}\cdot \mathop{\mathop{{\sum}^*}_{\chi\bmod
q}}_{\chi(-1)=-1}q^2\cdot (2\bar{\chi}(2)-4)\cdot L(m,\chi)L(1,\bar{\chi})\\
&&=\frac{q^{3-m}}{\phi(q)}\cdot \sum^{m+1}_{l=0}r_{m,1,l}\cdot q^{l-m}\cdot \frac{2^l-2^{m+2}-6}{2^{m-1}+1}
\prod_{p|q}\left(1-\frac{1}{p}\right)\left(1-\frac{1}{p^l}\right)\left(1-\frac{1}{p^{l-m}}\right)
\end{eqnarray*}
This completes the proof of Theorem 2.

In a similar way, we will deduce the identities involving generalized Hardy sums $s_2(\bar{a}b,m,n,q)$, $s_3(\bar{a}b,n,q)$, $s_5((\overline{2a-1})(2b-1),m,q)$ with Kloosterman sums respectively.

If $q$ is a square-full even number and  $m\equiv n\equiv 1(\bmod 2)$, we have
\begin{eqnarray*}
&&\mathop{{\sum}'}_{a=1}^{q}\mathop{{\sum}'}_{b=1}^{q}K(a,q)K(b,q)s_2(\bar{a}b,m,n,q)\\
&&=2^{m}\mathop{{\sum}'}_{a=1}^{q}\mathop{{\sum}'}_{b=1}^{q}K(a,q)K(b,q)S(2\bar{a}b,m,n,q)-\mathop{{\sum}'}_{a=1}^{q}\mathop{{\sum}'}_{b=1}^{q}K(a,q)K(b,q)S(\bar{a}b,m,n,q)\\
&&=-\frac{4m!n!}{(2\pi
i)^{m+n}q^{m+n-1}}\cdot \sum_{d|q}\frac{d^{m+n}}{\phi(d)}\cdot \mathop{\sum_{\chi\bmod
d}}_{\chi(-1)=-1}\mathop{{\sum}'}_{a=1}^{q}\left|\chi(a)K(a,q)\right|^2\cdot (2^m\bar{\chi}(2)-1)\cdot L(m,\chi)L(n,\bar{\chi})\\
&&=-\frac{4m!n!}{(2\pi
i)^{m+n}q^{m+n-1}}\cdot \frac{q^{m+n}}{\phi(q)}\cdot \mathop{\mathop{{\sum}^*}_{\chi\bmod
q}}_{\chi(-1)=-1}q^2\cdot (2^m\bar{\chi}(2)-1)\cdot L(m,\chi)L(n,\bar{\chi})\\
&&=\frac{q^{4-m-n}}{\phi(q)}\sum^{m+n}_{l=0}r_{m,n,l}\cdot (\frac{2^l-2^{m+n}-2^{n-m}-3}{2^{n-m}+1})
\cdot q^{l}\cdot \prod_{p|q}\left(1-\frac{1}{p}\right)\left(1-\frac{1}{p^l}\right)\left(1-\frac{1}{p^{l-m}}\right).
\end{eqnarray*}

If $q$ is a square-full odd number and  $n\equiv 1(\bmod 2)$, we get
\begin{eqnarray*}
&&\mathop{{\sum}'}_{a=1}^{q}\mathop{{\sum}'}_{b=1}^{q}K(a,q)K(b,q)s_3(\bar{a}b,n,q)\\
&&=2\mathop{{\sum}'}_{a=1}^{q}\mathop{{\sum}'}_{b=1}^{q}K(a,q)K(b,q)S(\bar{a}b,1,n,q)-4\mathop{{\sum}'}_{a=1}^{q}\mathop{{\sum}'}_{b=1}^{q}K(a,q)K(b,q)S(2\bar{a}b,1,n,q)\\
&&=-\frac{4n!}{(2\pi
i)^{1+n}q^{n}}\cdot \sum_{d|q}\frac{d^{n+1}}{\phi(d)}\cdot \mathop{\sum_{\chi\bmod
d}}_{\chi(-1)=-1}\mathop{{\sum}'}_{a=1}^{q}\left|\chi(a)K(a,q)\right|^2\cdot (2-4\bar{\chi}(2))\cdot L(1,\chi)L(n,\bar{\chi})\\
&&=-\frac{4n!}{(2\pi
i)^{1+n}q^{n}}\cdot \frac{q^{n+1}}{\phi(q)}\cdot \mathop{\mathop{{\sum}^*}_{\chi\bmod
q}}_{\chi(-1)=-1}q^2\cdot (2-4\bar{\chi}(2))\cdot L(1,\chi)L(n,\bar{\chi})\\
&&=\frac{q^{3-n}}{\phi(q)}\cdot \sum^{n+1}_{l=0}r_{1,n,l}\cdot q^{l}\cdot \left(\frac{5\cdot 2^n-2^{l+1}+6}{2^{n-1}+1}\right)
\cdot \prod_{p|q}\left(1-\frac{1}{p}\right)\left(1-\frac{1}{p^l}\right)\left(1-\frac{1}{p^{l-n}}\right).
\end{eqnarray*}

If $q$ is a square-full odd number and$n\equiv 1(\bmod 2)$, it is clear that if $a$ pass through a reduced residue system mod
 $q$, then $2a-1$ also pass through a reduced residue system mod $q$, that is
\begin{eqnarray*}
\mathop{{\sum}'}_{a=1}^{q}\chi(2a-1)K(2a-1,q)&=&\mathop{{\sum}'}_{a=1}^{q}\chi(a)K(a,q)=\tau^2(\chi).
\end{eqnarray*}
and
\begin{eqnarray*}
&&\mathop{\mathop{{\sum}'}_{a=1}^{q}}_{(2a-1,q)=1}\mathop{\mathop{{\sum}'}_{b=1}^{q}}_{(2b-1,q)=1}K(2a-1,q)K(2b-1,q)s_5((\overline{2a-1})(2b-1),m,q)\\
&&=2\cdot \mathop{\mathop{{\sum}'}_{a=1}^{q}}_{(2a-1, q)=1}\mathop{\mathop{{\sum}'}_{b=1}^{q}}_{(2b-1,q)=1}K(2a-1,q)K(2b-1,q)S((\overline{2a-1})(2b-1),m,q)\\
&&-4\cdot \mathop{\mathop{{\sum}'}_{a=1}^{q}}_{(2a-1,q)=1}\mathop{\mathop{{\sum}'}_{b=1}^{q}}_{(2b-1,q)=1}K(2a-1,q)K(2b-1,q)S(\bar{2}(\overline{2a-1})(2b-1),m,q)\\
&&=-\frac{4m!}{(2\pi i)^{m+1}q^{m}}\cdot \left[2\cdot
\sum_{d|q}\frac{d^{m+1}}{\phi(d)}\mathop{\sum_{\chi\bmod
d}}_{\chi(-1)=-1}\left|\mathop{{\sum}'}_{a=1}^{q}\chi(2a-1)K(2a-1,q)\right|^2\cdot L(m,\chi)L(1,\bar{\chi})\right.\\
&&\quad\left.-4\cdot\sum_{d|q}\frac{d^{m+1}}{\phi(d)}\mathop{\sum_{\chi\bmod
d}}_{\chi(-1)=-1}\left|\mathop{{\sum}'}_{a=1}^{q}\chi(2a-1)K(2a-1,q)\right|^2\cdot \chi(2)L(m,\chi)L(1,\bar{\chi})\right]\\
&&=-\frac{4m!}{(2\pi i)^{m+1}q^{m}}\cdot \frac{q^{m+1}}{\phi(q)}\cdot
\mathop{\mathop{{\sum}^*}_{\chi\bmod
q}}_{\chi(-1)=-1}q^2\cdot (2-4\chi(2))\cdot L(m,\chi)L(1,\bar{\chi})\\
&&=\frac{q^{3-m}}{\phi(q)}\cdot\sum^{m+1}_{l=0}r_{m,1,l}\cdot \left(\frac{5\cdot 2^m-2^{l+1}+6}{2^{m-1}+1}\right)\cdot
q^{l}\cdot \prod_{p|q}\left(1-\frac{1}{p}\right)\left(1-\frac{1}{p^l}\right)\left(1-\frac{1}{p^{l-m}}\right)
\end{eqnarray*}

\bigskip

This completes the proof of Theorem.

\bigskip

\bigskip

\section*{References}

[1] Apostol, T. M.: Modular Functions and Dirichlet Series in Number
Theory, Springer-Verlag, New York, 1976.

[2] Carlitz, L.: The reciprocity theorem of Dedekind sums. Pacific J.
Math., 3, 513-522 (1953).

[3] B.C. Berndt, Analytic Eisenstein series, theta-functions, and
series relations in the spirit of Ramanujan, J. Reine Angew. Math.
303/304 (1978) 332-365.

[4] R. Sitaramachandrarao, Dedekind and Hardy sums, Acta Arith. 48
(1978) 325–340.

[5] B.C. Berndt: Analytic Eisenstein series, theta-functions, and
series relations in the spirit of Ramanujan. J. Reine Angew. Math.
303/304, 332-365 (1978)

[6] B.C. Berndt and L.A.:Analytic Properties of Arithmetic Sums
Arising in the Theory of the Classical Theta-Functions, SIAM Journ
Math. Anal.1984,15(1):143-150.

[7] R. Sitaramachandrarao, Dedekind and Hardy sums, Acta Arith. 48
(1978) 325-340.

[8] M.R. Pettet, R. Sitaramachandrarao, Three-term relations for
Hardy sums, J. Number Theory 25 (1987) 328-339.

[9] Huaning Liu, Wenpeng Zhang,.Generalized Dedekind sums and
Generalized  Hardy sums, Acta Mathematica Sinica, Chinese Series
2006, 49(5): 999-1008.

[10] Y. Liu, W. Zhang: A hybrid mean value related to the Dedekind
sums and Kloosterman sums. Sci. China, Math. 53 (2010), 2543–2550

[11] H. Zhang, W.P. Zhang, On the identity involving certain Hardy
sums and Kloosterman sums, J.Inequal. Appl. 2014 (2014) 52.

[12] Wen Peng, Tianping Zhang. Some identities involving certain
Hardy sum and Kloosterman sum, Journal of Number
Theory,165(2016)355-362

[13] Tom M. Apostol, Introduction to Analytic Number Theory,
Springer-Verlag, New York, 1976.

\end{document}